\documentclass[10pt,leqno]{amsart}
\usepackage{amscd}
\usepackage{amssymb}
\usepackage{amsfonts}
\usepackage{latexsym}
\usepackage{verbatim}
\usepackage{amsthm}
\numberwithin{equation}{section}

\theoremstyle{plain}
\newtheorem{theorem}{Theorem}[section]
\newtheorem{definition}[theorem]{Definition}
\newtheorem{lemma}[theorem]{Lemma}
\newtheorem{prop}[theorem]{Proposition}
\newtheorem{cor}[theorem]{Corollary}
\newtheorem{rem}[theorem]{Remark}

\newcommand{\w}[1]{\widetilde{#1}}

\renewcommand{\o}{\omega}
\renewcommand{\a}{\alpha}
\renewcommand{\w}{\wedge}
\begin{document}
\title{Torsion of SU(2)-structures and Ricci curvature in dimension 5}
\author{Lucio Bedulli and Luigi Vezzoni}
\date{\today}
\address{Dipartimento di Matematica \\ Universit\`a di Torino\\
Via Carlo Alberto 10 \\
10123 Torino\\ Italy} \email{luigi.vezzoni@unito.it}
\address{Dipartimento di Matematica e Applicazioni per l'Architettura - Universit\`a di Firenze\\
Piazza Ghiberti 27\\50122 Firenze\\Italy}
\email{bedulli@math.unifi.it} 
\thanks{This work was supported by the Project M.I.U.R. ``Riemannian metrics and differentiable manifolds'' and by G.N.S.A.G.A.
of I.N.d.A.M.}
\begin{abstract}
Following the approach of Bryant \cite{B2}, we study the intrinsic
torsion of an SU$(2)$-structure on a $5$-dimensional manifold deriving an explicit expression for the
Ricci and the scalar curvature in terms of torsion forms and its derivative.
As a consequence of this formula  we prove that the $\alpha$-Einstein condition forces some special
SU(2)-structures to be Sasaki-Einstein.
\end{abstract}
\maketitle
\newcommand\C{{\mathbb C}}
\newcommand\R{{\mathbb R}}
\newcommand\Z{{\mathbb Z}}
\newcommand\T{{\mathbb T}}
\newcommand\GL{{\rm GL}}
\newcommand\SL{{\rm SL}}
\newcommand\SO{{\rm SO}}
\newcommand\Sp{{\rm Sp}}
\newcommand\U{{\rm U}}
\newcommand\SU{{\rm SU}}
\newcommand{\Gdue}{{\rm G}_2}
\newcommand\re{\,{\rm Re}\,}
\newcommand\im{\,{\rm Im}\,}
\newcommand\id{\,{\rm id}\,}
\newcommand\tr{\,{\rm tr}\,}
\renewcommand\span{\,{\rm span}\,}
\newcommand\Ann{\,{\rm Ann}\,}
\newcommand\Hol{{\rm Hol}}
\newcommand\Ric{{\rm Ric}}

\newcommand\nc{\widetilde{\nabla}}
\renewcommand\d{{\partial}}
\newcommand\dbar{{\bar{\partial}}}
\newcommand\s{{\sigma}}
\newcommand\sd{{\bigstar_2}}
\newcommand\K{\mathbb{K}}
\renewcommand\P{\mathbb{P}}
\newcommand\D{\mathbb{D}}

\newcommand\f{{\varphi}}
\newcommand\g{{\frak{g}}}
\renewcommand\k{{\kappa}}
\renewcommand\l{{\lambda}}
\newcommand\m{{\mu}}
\renewcommand\O{{\Omega}}
\renewcommand\t{{\theta}}
\newcommand\ebar{{\bar{\varepsilon}}}
\section*{Introduction}
In 1960 Sasaki introduced in \cite{Sas} a new class of
contact-metric structures which can be considered as an
odd-dimensional counter-part of K\"ahler structures. This kind of
geometry became known as \emph{Sasakian geometry} and it is
present today in many mathematical and physical contexts. In
Sasakian geometry Einstein metrics play a central role and
Sasaki-Einstein manifolds arise in many physical models. As
general references for these topics see e.g. \cite{BL}, \cite{BG},
\cite{BGM}, \cite{BGN}, \cite{FOW}, \cite{MSY}, \cite{MSY1} and
the references therein.\\ Since in dimension $5$ Sasakian-Einstein
metrics correspond to Killing spinors (see \cite{FK}), it is
rather natural to study the larger class of SU(2)-structures
induced by {\it generalized} Killing spinors. These structures
were firstly investigated and called \emph{Hypo-structures} by
Conti and Salamon in \cite{CS}, where they prove that any analytic
Hypo-manifold can be realized as a hypersurface of a Calabi-Yau
threefold.\\ In terms of differential forms a Hypo-structure is
determined by a nowhere vanishing $1$-form $\alpha$ and a triple
of $2$-forms $(\o_1,\o_2,\o_3)$ satisfying $$
\o_i\w\o_j=\delta_{ij\,}v\,\quad d\o_1=0\,,\quad
d(\o_2\w\a)=0\,,\quad d(\o_3\w\a)=0\,. $$ where $v$ is a $4$-form
such that $v\w\a\neq 0$ everywhere.\\

In \cite{FIMU} the authors introduce two new types of
SU(2)-structures on $5$-manifolds: {\it nearly-Hypo} structures
are the natural structures inherited by an hypersurface of a
nearly K\"ahler SU(3)-manifold, while {\it double-Hypo} structures
are nearly-Hypo and Hypo simultaneously.

In this paper, following the same approach used by Bryant in
\cite{B2} to compute the Ricci tensor of a G$_2$-structure, we
write down an explicit formula for the scalar curvature and the
Ricci tensor of the metric induced by an SU(2)-structure on a
$5$-manifold in terms of the intrinsic torsion (Theorems
\ref{scal} and \ref{Ricci}). As a direct consequence of the
formula of the scalar curvature, we have that if the Ricci tensor
of a Hypo-structure $(\alpha,\omega_1,\omega_2,\omega_3)$
satisfies $$ \mbox{Ric}(R_{\alpha},R_{\alpha})=4\,, $$ where
$R_{\alpha}$ is the Reeb vector field of $\alpha$, then the
Hypo-structure is Sasaki $\alpha$-Einstein. This result slightly strengthen a
previous result by Conti and Salamon (see \cite{CS}).
\\
The formula for the Ricci tensor has as a direct application the study of $\alpha$-Einstein metrics on contact-Hypo manifold.
The $\alpha$-Einstein metrics were introduced by Okumura in \cite{Oku} in the context of contact-metric geometry and they
are characterized by the equation
\[{\rm Ric} = \mu g+ \lambda \alpha \otimes \alpha\,,
\]
where $\lambda$ and $\mu$ are constant.
Sasaki $\alpha$-Einstein metrics seem to be a natural generalization of K\"ahler-Einstein metrics to the odd dimension
(see e.g \cite{BGM}).\\
We prove that the $\alpha$-Einstein condition forces a double-Hypo structure to be Sasaki-Einstein
(Proposition \ref{double}). Finally, as a corollary, we prove that if the almost K\"ahler cone of a 5-dimensional $\alpha$-Einstein 
SU(2)-manifold inherits a symplectic half-flat structure (see \cite{dBT}, \cite{CT} and \cite{BV}), then it is a 
Sasaki $\alpha$-Einstein manifold (Corollary \ref{SGCY}).

The present paper is organized as follows: In section 1 we recall
some basic facts on SU(2)-structures and set up the algebraic
preliminaries needed in the sequel. In section 2 we recall the
properties of the intrinsic torsion of an SU(2)-structure proving
some new formulae which will be useful in the next part of the
paper. Section 3 is devoted to the main result. We describe the
computational steps needed to reach it (and carried out with the
aid of {\sc Maple}) and we write down the formulae for the scalar
curvature and the Ricci tensor. Then we prove the consequences
obtained imposing the $\alpha$-Einstein condition. \\
{\sc Acknowledgments.} The authors are grateful to Robert Bryant for
supplying them with the computer programs he used to perform the symbolic
computations in the G$_2$-case and to Diego Conti for some observations which helped to
improve the presentation of this paper. Finally the authors are grateful to an anonimous referee
for pointing out some mistakes.
\\
\newline
{\sc Notation.} Given a manifold $M$, we denote by $\Lambda^rM$ the
space of smooth $r$-forms on $M$.\\
When a coframe $\{e^1,\dots,e^n\}$ is given, we
will denote the $r$-form $e^{i_{1}}\wedge\dots\wedge e^{i_{r}}$ by
$e^{i_{1}\dots i_{r}}$.\\
Furthermore when a contact form $\alpha$ is fixed, we will denote by $\beta^T$ the projection
of an arbitrary differential form $\beta$ onto the contact distribution $\xi=\ker\alpha$.\\
The symbol $\langle\cdot\,,\cdot\rangle$ will denote the scalar product induced on exterior forms
by a Riemannian metric.\\
Finally in the indicial expression the symbol of sum over repeated indices
is omitted.

\section{Five-dimensional SU$(2)$-structures}
Let $M$ be a 5-dimensional smooth manifold and let
$\mathcal{L}(M)\to M$ be the $\mbox{GL}(5)$-bundle of linear
frames on $M$. An $\mbox{SU}(2)$-structure on $M$ is by definition an
$\mbox{SU}(2)$-reduction of $\mathcal{L}(M)$. In terms of
differential forms an $\mbox{SU}(2)$-structure may be characterized
as follows
\begin{prop}[\cite{CS}]
$\mbox{\emph{SU}}(2)$-structures on $M$ are in one-to-one
correspondence with quadruples
$(\alpha,\omega_1,\omega_2,\omega_3)$, where $\alpha$ is a nowhere
vanishing $1$-form, $\omega_1,\omega_2,\omega_3\in\Lambda^2M$
satisfy
$$
\omega_i\wedge\omega_j=\delta_{ij}\,v\quad \mbox{for }i=1,2,3\,
$$
for some $4$-form $v$ with $v\wedge \alpha\neq 0$, and
\[
\iota_X \omega_1=\iota_Y\omega_2 \Longrightarrow \omega_3(X,Y)\geq
0\,.
\]
\end{prop}
An $\mbox{SU}(2)$-structure $(\alpha,\omega_1,\omega_2,\omega_3)$ on
$M$ singles out a rank 4 distribution $\xi=\ker\alpha \subset TM$.
Note that for any $r=1,2,3$, the pair $(\xi,\omega_r)$ is a
symplectic bundle over $M$. Furthermore there exists a unique vector
field $R_{\alpha}$ on $M$ satisfying
$$
\alpha(R_{\alpha})=1\,,\quad \iota_{R_{\alpha}}\omega_1=0\,.
$$
In analogy with the terminology used in contact geometry, we will
refer to $R_{\alpha}$ as the \emph{Reeb vector field} associated to
$(\alpha,\omega_1,\omega_2,\omega_3)$. Note that from the definition
we also have
$$
\iota_{R_{\alpha}}\omega_2=\iota_{R_{\alpha}}\omega_3=0\,.
$$
\begin{definition}
{\em A differential form $\gamma$ on $M$ is said to be
{\em $\alpha$-transversal} if it satisfies $\iota_{R_{\alpha}}\gamma=0$.
The set of $\alpha$-transversal $p$-forms on $M$ is denoted by
$\Lambda^p_0M$. Analogously $S^p_0(M)$ will denote the set of
$\alpha$-transversal symmetric $p$-tensors defined in the same way.
}
\end{definition}
\begin{rem}
\emph{ If we identify the vector bundle $\xi^*$ dual to $\xi$ with
the subbundle of $T^*M$ whose fibre over $x$ is $\{\phi\in
T_x^*M\,\,|\,\,\phi(R_{\alpha})=0\}$, then $\Lambda^p_0M$ is
identified with $\Gamma(\Lambda^{p}\xi^*)$.}
\end{rem}
We define the operators $\bigstar_r$, $r=1,2,3$ on the transversal
forms
$$
\bigstar_{r}\colon\Lambda^{j}_0M\to\Lambda_{0}^{4-j}M
$$
by means of the relations
$$
\gamma\wedge\bigstar_r\,\beta=\omega_r(\gamma,\beta)\frac{\omega_r^2}{2}\,,
$$
for $r=1,2,3$, where the $\omega_r$'s are extended to exterior forms
in the usual way.
\begin{lemma}
Let
$$
J_r\colon \Lambda^1M\to \Lambda^1M\,,\mbox{ for }r=1,2,3\,,
$$
be the $C^{\infty}(M)$-linear endomorphisms defined by
$$
\begin{aligned}
& J_{1}(\phi)=\bigstar_1(\omega_3\wedge\bigstar_1(\omega_2\wedge\phi))\,,\\
& J_{2}(\phi)=\bigstar_2(\omega_1\wedge\bigstar_2(\omega_3\wedge\phi))\,,\\
&
J_{3}(\phi)=\bigstar_3(\omega_2\wedge\bigstar_3(\omega_1\wedge\phi))\,,
\end{aligned}
$$
for any $\phi\in \Lambda^1_0M$ and by
$$
J_{1}(\alpha)=J_{2}(\alpha)=J_{3}(\alpha)=0\,.
$$
Then for $r=1,2,3$ one has
\begin{itemize}
\item $J_r^2=-I+ R_{\alpha}\otimes\alpha\,$;
\vskip0.1cm
\item $\omega_r(J_r\beta,\gamma)=-\omega_r(\beta,J_r\gamma)\,$ for every $1$-form $\beta, \gamma$.
\end{itemize}
\end{lemma}
\begin{proof}
The statement is a consequence of the real version of Schur's lemma
(the proof is analogous to that of Proposition 2.1 of \cite{BV}).
\end{proof}

Every $J_{r}$ induces an endomorphism of $TM$ (we denote it with the
same symbol) in the following way
\begin{enumerate}
\item[1.] if $X$ is a smooth section of $\xi$, then we set
$J_{r}(X):=-\sharp_r^{-1}J_r(\sharp_{r}X)$, where $\sharp_r\colon \xi\to\xi^*$ is
the duality on $\xi$ induced by $\omega_r$, \vskip0.1cm
\item[2.] if $X=R_{\alpha}$ we set $J_r(X)=0$.
\end{enumerate}
In this way each $J_r$ is an $\omega_r$-compatible bundle complex structure on $\xi$.\\
Note that from the definition one easily obtains the 
quaternionic identities satisfied by $J_r \in {\rm End}(TM)$:
\[
J_rJ_s=-J_sJ_r\,,\quad \mbox{ for }r,s=1,2,3\,,\,r\neq s
\]
and
\[
J_1J_2=J_3\,.
\]
At the dual level the $J_r$'s anticommute, but the composition satisfies $J_1 J_2 = -J_3$.

Furthermore we fix on $M$ the Riemannian metric $g$ defined by
$$
g=g^T+\alpha\otimes\alpha\,,
$$
where
$$
g^T(X,Y)=\omega_1(X,J_1Y)=\omega_2(X,J_2Y)=\omega_3(X,J_3Y)\,.
$$
Note that for any $X,Y\in\Gamma(\xi)$ we have
$$
g^T(J_1X,J_1Y)=g^T(J_2X,J_2Y)=g^T(J_3X,J_3Y)=g^T(X,Y)\,.
$$
Another direct consequence is that
$$
g(J_rX,J_rY)=g(X,Y)-\alpha(X)\alpha(Y)\, \mbox{ for } r=1,2,3\,.
$$
The metric $g$ together with the orientation defined by
$\alpha\wedge\omega_1^2$ induces the Hodge star operator in the
usual way. Finally we denote by $*^T$ the transverse Hodge operator
acting on the transverse forms so that
$$
\eta\wedge *^T\nu=g^T(\eta,\nu)\,\frac{\omega_1^2}{2}\,.
$$
Note that
$$
*^T\omega_r=\omega_r\,\,\mbox{ for }r=1,2,3
$$
and that for any transverse $p$-form $\gamma$ we have
$$
*^T\gamma=*(\alpha\w\gamma)\,.
$$
\subsection{The standard model}
Let $e^1,\ldots,e^5$ be the coframe dual to the canonical basis of
$\R^5$. Then
\begin{equation}
\label{standard}
\begin{aligned}
{\mathbf \alpha} & =  e^5\,,\\
\omega_1 & =e^{12}+e^{34}\,,\\
\omega_2 & =e^{13}-e^{24}\,,\\
\omega_3 & =e^{14}+e^{23}\,.
\end{aligned}
\end{equation}
define a linear SU(2)-structure on $\R^5$. In fact, given any linear
SU(2)-structure on a vector space $V$, we can find a basis of $V^*$
with respect to which the structure forms take the standard form
\eqref{standard} (see \cite{CS}). Therefore it is useful to
introduce the following notation:
\[
\omega_r=\frac{1}{2}\epsilon^r_{ij}\,e^i\wedge e^j\,.
\]
The endomorphisms $J_1,J_2,J_3$ induced by the standard structure
act on the canonical basis $e_1,\dots,e_5$ as follows
\begin{eqnarray*}
J_1(e_1)=e_2  & J_2(e_1)=e_3 & J_3(e_1)=e_4\\
J_1(e_3)=e_4  & J_2(e_4)=e_2 & J_3(e_2)=e_3.
\end{eqnarray*}
Using this standard model one can easily check that, given an
SU(2)-structure $(\alpha,\omega_1,\omega_2,\omega_3)$ on a
$5$-dimensional manifold $M$,
\begin{equation}
\label{formuletta} \bigstar_{r}\phi=\phi\wedge\omega_r\,,
\end{equation}
for any $r=1,2,3$ and transverse $1$-form $\phi$ on $M$.

\subsection{Decomposition of the Lie algebra $\mathfrak{so}(5)$}
We use the $\epsilon$-notation introduced above to obtain the
decomposition of the Lie algebra $\mathfrak{so}(5)$ of
skew-symmetric $5\times 5$  matrices in irreducible SU(2)-modules.
Indeed
\begin{equation}
\label{2-tensors} \mathfrak{so}(5)\simeq
\mathfrak{su}(2)\oplus[\![\R^4]\!]\oplus[\R]_1\oplus[\R]_2\oplus[\R]_3\,,
\end{equation}
where a matrix $A=(a_{ij})$ lies in $\mathfrak{su}(2)$ if and only
if
$$
\begin{cases}
\epsilon^r_{ij}a_{ij}=0\\
a_{i5}=a_{5i}=0\,;
\end{cases}
$$
for every $v=(v_1,v_2,v_3,v_4)\in\R^4$

$$
[\![v]\!]=
\begin{pmatrix}
&0 &0 &0 &0 &v_1\\
&0 &0 &0 &0 &v_2\\
&0 &0 &0 &0 &v_3\\
&0 &0 &0 &0 &v_4\\
&-v_1 &-v_2 &-v_3 &-v_4 &0\\
\end{pmatrix}
$$
and for any $t\in\R$
$$
([\,t\,]_r)_{ij}=t\epsilon^r_{ij}\,.
$$
Note that we can alternatively write in compact form
$$
[\![v]\!]_{ij}=\eta_{ijk} v_{k}\,,
$$
using the $\eta$-symbol
\begin{equation}
\label{eta}
\eta_{ijk}=\delta_{ik}\delta_{j5}-\delta_{jk}\delta_{i5}\;\mbox{ for
}i,j=1\dots,5\,,\;k=1,\dots,4
\end{equation}
we will need later.
\section{Intrinsic torsion and special $\mbox{SU}(2)$-structures}
Since the natural action of SU(2) on $\Lambda^p(\R^5)^*$ for every $p$, once an SU(2)-structure on a 5-manifold $M$ is fixed, we have
a natural splitting of the space of forms of each
degree. More precisely we
have the following decomposition in irreducible SU(2)-modules:
$$
\begin{aligned}
&\Lambda^{1}M=\langle\alpha\rangle\oplus\Lambda^{1}_0M\,,\\
&\Lambda^2M=\alpha\wedge
\Lambda^1_0M\oplus_{r=1}^3\langle\,\omega_r\rangle\oplus
\Lambda^{2}_{3}M\,,\\
&\Lambda^{3}M=\Lambda^{3}_0M\oplus_{r=1}^3\langle\,\alpha\wedge\omega_r\rangle\oplus
\alpha\wedge\Lambda^{2}_{3}M\,,
\end{aligned}
$$
where
$$
\Lambda^{2}_3M=\{\sigma\in\Lambda^2_0M\,\,|\,\,\sigma\wedge\omega_r=0\,\,\,\mbox{for}\,\,r=1,2,3\}\,.
$$
The previous decomposition allows us to define also a projection
$$
\mbox{E}\colon \Lambda^2M\to \Lambda^{2}_3M
$$
explicitly defined by
\begin{equation}
\label{E}
\mbox{E}(\phi)=\phi^T-\sum_{r=1}^3\frac{1}{2}*(\phi^T\wedge\omega_r\wedge\alpha)\,\omega_r\,,
\end{equation}
where $\phi^T$ denotes the projection of $\phi$ onto
$\Lambda^{2}_0M$, i.e.
$$
\phi^T:=\phi-\alpha\wedge\iota_{R_{\alpha}}\phi\,.
$$
\begin{rem}\emph{
Since E is the projection on the -1 eigenspace of the diagonalizable operator $*^T$,
the operator E restricted to $\Lambda_0^2M$ and $*^T$ commute, i.e.
\[
{\rm E}(*^T\beta)=*^T{\rm E}(\beta)
\]
for every $\beta \in \Lambda^2_0 M$.
Moreover, if $\psi$ is an arbitrary 3-form on $M$, then we
immediately have
\begin{equation}
\label{prop*} \mbox{E}(*\psi)=
*^T\mbox{E}(\iota_{R_{\alpha}}\psi)\,.
\end{equation}
}
\end{rem}
\begin{rem}\label{aggiunto}
\emph{ Note that the elements of $\Lambda^2_3M$ are the sections of
a subbundle of $\Lambda^2(T^*M)$ isomorphic to the bundle associated
to the SU(2)-reduction $\mathcal{Q}$ with respect to the adjoint
representation of SU(2).}
\end{rem}
In the sequel we will use the following
\begin{prop}
Let $\sigma\in \Lambda^{2}_3M$, then
\begin{enumerate}
\item[1.] $*\sigma=-\sigma\wedge\alpha$,
\vskip0.1cm
\item[2.] $J_r(\sigma)=\sigma$ for $r=1,2,3$.
\end{enumerate}
\end{prop}
\begin{proof}
Since any element of the Lie algebra $\mathfrak{su}(2)$ is
SU(2)-conjugated to an element of a fixed Cartan subalgebra, remark
\ref{aggiunto}  implies that for any $x\in M$  there exists an
SU(2)-local frame $e_1,\dots,e_5$ near $x$, such that
$$
\sigma=e^{12}-e^{34}
$$
and the claim follows.
\end{proof}
According to the decomposition of the exterior algebra the
derivatives of the structure forms split as
$$
\begin{aligned}
d\omega_r&=\nu_r\wedge\omega_r+\sum_{j=1}^3f_{rj}\alpha\wedge\omega_j+\alpha\wedge\sigma_r\,,\\
d\alpha&=\alpha\wedge\nu_4+\sum_{i=1}^3 \phi_i\,\omega_i+\sigma_4\,,
\end{aligned}
$$
where $\nu_i\in\Lambda^1_0M$, $\sigma_{i}\in\Lambda^2_3M$, for
$i=1,\dots,4$ and $\phi_i$, $f_{ij}$ are smooth functions. Imposing
$d^2=0$ one has
$$
\begin{aligned}
&f_{11}=f_{22}=f_{33}\,,\\
&f_{ij}=-f_{ji}\mbox{ for }i\neq j\,.
\end{aligned}
$$
We will refer to  $\{\nu_i,\sigma_j,\phi_r,f_{uv}\}$ as the
\emph{torsion forms} of the SU(2)-structure.
\subsection{Decomposition of symmetric 2-tensors}
In order to write the Ricci tensor of a 5-dimensional SU(2)-manifold
in terms of its torsion forms, we must decompose the space of
symmetric 2-tensors on $M$ in irreducible SU(2)-modules. We have
\begin{equation}\label{decsym}
S^2(M)=\langle\, g^T\rangle\oplus\langle\alpha\otimes\alpha\rangle
\oplus_{i=1}^3\Sigma_i(M)\oplus(\alpha\odot\Lambda^1_0M)\,.
\end{equation}
where
$$
\begin{aligned}
&\Sigma_{1}(M)=\{h\in S^2_0(M)\,\,|\,\,J_1(h)=h,\,J_2(h)=J_3(h)=-h\}\,,\\
&\Sigma_{2}(M)=\{h\in S^2_0(M)\,\,|\,\,J_2(h)=h,\,J_1(h)=J_3(h)=-h\}\,,\\
&\Sigma_{3}(M)=\{h\in
S^2_0(M)\,\,|\,\,J_3(h)=h,\,J_1(h)=J_2(h)=-h\}\,.
\end{aligned}
$$
Let
$$
\iota_{r}\colon \Sigma_r(M)\to\Lambda^2_{3}M
$$
be defined by
\begin{equation}
\label{iota} \iota_{r}(h_{lm}\,e^{l}\!\otimes\!
e^m)=\frac{1}{2}\epsilon^r_{ik}h_{kj}e^i\wedge e^j\,.
\end{equation}
It is immediate to verify that every $\iota_r$ is an isomorphism of SU(2)-representations.
\subsection{The almost K\"ahler cone and special SU(2)-structures }
\label{akcone} In order to consider some interesting kind of
SU(2)-structure on 5-manifolds, we first take the more general point
of view of U($n$)-structures on $(2n+1)$-manifolds. A U($n$)
structure on a $(2n+1)$-dimensional manifold $M$ is determined by a
triple $(\alpha,J,\omega)$, where $\alpha$ is a nowhere vanishing
1-form on $M$, $\omega$ is a 2-form such that
\[
\alpha\wedge\omega^n\neq 0\,,
\]
and $J\in \rm{End}(TM)$ is such that
\[
J^2=-I+\alpha \otimes R_\alpha\,,
\]
where $R_\alpha$ is the {\em Reeb vector field} (i.e.
$\alpha(R_\alpha)=1$ and $\iota_{R_\alpha} \omega=0$). Any
U($n$)-structure on $M$ induces a U($n+1$)-structure on the cone
$C(M)=M \times \R^+_t$ specified by
\begin{equation}
\label{kappa} \kappa=t^2\omega+t\alpha \wedge dt
\end{equation}
and the $\kappa$-compatible almost complex structure $\widetilde{J}$ defined by
\[
\widetilde{J}X=
\begin{cases}
JX \quad \;\,\,
\mbox{if } X\in \Gamma(\ker \alpha)\\
-t\frac{\partial}{\partial t} \quad \mbox{if } X=R_\alpha\,.
\end{cases}
\]
Note that the 2-form $\kappa$ is closed (and hence symplectic) if
and only if $\alpha$ and $\omega$ are related by
\[
d\alpha=-2\omega\,.
\]
In this case $\alpha$ is a contact form on $M$ and $(\kappa,
\widetilde{J})$ is an {\em almost-K\"ahler} structure on $C(M)$. A
U($n$)-structure is said to be {\em Sasakian $($Sasaki-Einstein$)$} if
$(C(M),\kappa,\widetilde{J})$ is a
K\"ahler (Calabi-Yau) manifold.

Now let us come back to the case of SU(2)-structures. First remark
that an SU(2)-structure $(\alpha,\omega_1,\omega_2,\omega_3)$ on a
5-dimensional manifold $M$ induces an SU(3)-structure on the cone
$C(M)$. In fact, once a U(3)-structure $(\kappa,J)$ on a
6-dimensional manifold $N$ is given, in order to specify an
SU(3)-structure it is sufficient to give a complex volume form
$\varepsilon\in\Lambda^{3,0}_JN$ satisfying
$$
\varepsilon\wedge\overline{\varepsilon}=-i\frac{4}{3}\,\kappa^3\,.
$$
Hence on the cone $C(M)$ we may take
$$
\varepsilon=t^2(\omega_2+i\omega_3)\wedge(t\alpha+idt)\,.
$$
This SU(3)-structure is integrable if and only if
$$
\begin{cases}
\begin{aligned}
&d\alpha=-2\omega_1\\
&d\omega_2=3\alpha\wedge\omega_3\\
&d\omega_3=-3\alpha\wedge\omega_2\,,
\end{aligned}
\end{cases}
$$
see e.g. \cite{CS}.

Here we list some special SU(2)-structures which have been studied
in the last years.

\begin{itemize}
\item \textbf{Hypo manifolds}: An $\mbox{SU}(2)$-structure $(\alpha,\omega_1,\omega_2,\omega_3)$
is said to be a Hypo-structure if the structure forms satisfy
$$
d\omega_1=0\,,\quad d(\alpha\wedge\omega_2)=0\,,\quad
d(\alpha\wedge\omega_3)=0\,.
$$
In terms of intrinsic torsion the Hypo condition reads as
$$
\begin{aligned}
&\nu_1=0\,,\quad f_{1j}=0\,,\,\,i=1,2,3\,,\quad \sigma_1=0\,,\\
&\nu_2=\nu_3=\nu_4\,,\quad \phi_2=\phi_3=0
\end{aligned}
$$
and the other torsion forms are arbitrary. Hypo structures were
first investigated by Conti and Salamon in \cite{CS}. The name is
due to the fact that a real hypersurface of a Calabi-Yau 3-fold
inherits a Hypo structure.
\vskip0.3cm
\item \textbf{Contact-Hypo manifolds}: A Hypo structure is called contact Hypo if further the 1-form $\alpha$ is a contact
form so that the SU(3)-structure on the cone $C(M)$ is actually
almost K\"ahler, i.e.
$$
d\alpha=-2\omega_1\,,\quad d(\alpha\wedge\omega_2)=0\,,\quad
d(\alpha\wedge\omega_3)=0\,.
$$
This special SU(2)-structures are the subject of the recent paper \cite{AFFU}.
In terms of torsion forms we have
\begin{equation}
\label{CHtorsion}
\begin{aligned}
&\nu_i=0\,,\,\,i=1,2,3,4\,,\quad f_{1j}=0\,,\,\,i=1,2,3\,,\\
&\sigma_1=\sigma_4=0\,,\quad\phi_1=-2\,,\,\phi_2=\phi_3=0\,.
\end{aligned}
\end{equation}
\vskip0.3cm
\item \textbf{Nearly Hypo manifolds}: These manifolds have been
introduced in \cite{FIMU}. In this case the structure equations are
$$
d\omega_2=3\alpha\wedge\omega_3\,,\quad
d(\alpha\wedge\omega_1)=-2\omega_1^2
$$
which in terms of torsion  forms are
$$
\begin{aligned}
&\phi_1=-2\,,\quad \phi_3=0\,,\quad \sigma_2=0\,,\\
&\nu_4=\nu_3=\nu_1\,,\quad \nu_2=0\,,\quad f_{23}=3\,,\quad f_{12}=f_{11}=0\,,
\end{aligned}
$$
the remaining torsion forms being arbitrary. Such a structure is
inherited by any hypersurface of a nearly-K\"ahler SU(3)-manifold.
\vskip0.3cm
\item \textbf{Double Hypo manifolds}: These manifolds have been
introduced in \cite{FIMU}, too. A double Hypo structure is an
SU(2)-structure which is both Hypo and nearly Hypo. This kind of
structures  are characterized by the following equations
$$
d\omega_1=0\,,\quad d(\alpha\wedge\omega_2)=0\,,\quad
d(\alpha\wedge\omega_1)=-2\omega_1^2\,,
$$
$$
d\omega_2=3\alpha\wedge\omega_3\,.
$$
In this case the only non-vanishing torsion forms are
$\phi_1,f_{23},\sigma_3,\sigma_4$, where
$$
\phi_1=-2\,,\quad f_{23}=3\,,
$$
and $\sigma_3$, $\sigma_4$ are arbitrary.
\vskip0.3cm
\item \textbf{Sasaki-Einstein manifolds}: A Sasakian manifold is
said to be Sasaki-Einstein if the induced Riemann metric is
Einstein. In dimension $5$ a Sasakian structure induced an
SU(2)-structure $(\alpha,\omega_1,\omega_2,\omega_3)$ satisfying
$$
d\alpha=-2\omega_1\,,\quad d\omega_2=3\alpha\wedge\omega_3\,,\quad
d\omega_3=-3\alpha\wedge\omega_2\,.
$$
In terms of torsion forms these conditions read as
$$
\phi_{1}=-2\,,\quad f_{23}=3
$$
and the other torsion forms vanish.
\end{itemize}
\subsection{Sympletic half-flat structures.}
Let $N$ be a 6-dimensional manifold. Any U(3)-structure $(\kappa,J)$
on $N$ induces a natural connection $\widetilde{\nabla}$, called the
Hermitian connection, uniquely determined by the following equations
$$
\widetilde{\nabla} J=0\,,\quad \widetilde{\nabla}\kappa=0\,,\quad
(T^{\widetilde{\nabla}})^{1,1}=0
$$
where $(T^{\widetilde{\nabla}})^{1,1}$ is the $(1,1)$-part of the
torsion of $\widetilde{\nabla}$. It turns out that the holonomy
group of this connection is contained in SU(3) if and only if there
exists $\varepsilon\in\Lambda^{3,0}_JN$ satisfying
$$
\begin{cases}
\varepsilon\wedge\overline{\varepsilon}=-i\frac{4}{3}\,\kappa^3\\
\overline{\partial}_J\varepsilon=0
\end{cases}
$$
(see e.g. \cite{dBT}). In this case we call $(N,\kappa,J,\varepsilon)$ 
a \emph{symplectic Calabi-Yau manifold} (\cite{dBT},
\cite{BV} \footnote{In \cite{dBT} and \cite{BV} such structures are named 
{\em Generalized Calabi-Yau}, but this terminology is misleading because 
it is widely used with a different meaning, see \cite{Hi}.}).                                                           
Requiring further that the
real part of $\varepsilon$ is closed, we obtain an interesting
subclass of manifolds lying in the intersection between
symplectic and half-flat geometry, indeed they are 
called {\em symplectic half-flat} manifolds in \cite{CT}.\\

Let us consider now a $5$-dimensional SU($2$)-manifold
$(M,\alpha,\omega_1,\omega_2,\omega_3)$. We have the following
\begin{lemma}
\label{SGCYcone}
Let $(\kappa,\widetilde{J},\varepsilon)$ be the
$\mbox{\emph{SU}}(3)$-structure on the cone $C(M)$ associated to
$(\alpha,\omega_1,\omega_2,\omega_3)$. Then $(\kappa,\widetilde{J},\varepsilon)$ is symplectic half-flat if and
only if $(\alpha,\omega_1,\omega_2,\omega_3)$ is contact-Hypo with
$f_{23}=3$ and $\sigma_3=0$.
\end{lemma}
\begin{proof}
As already observed, the 2-form $\kappa$ defined by \eqref{kappa} is
closed if and only if $d\alpha=-2\,\omega_1$. This implies
$d\omega_1=0$. In terms of torsion forms:
\[
\begin{array}{ccc}
\phi_1=-2&\phi_2=\phi_3=0& \\
f_{1r}=0& \nu_1=\nu_4=0&\sigma_1=\sigma_4=0\,.
\end{array}
\]
Now
\[
\begin{aligned}
d\,\mbox{Re}\,\varepsilon=&-3t^2\o_2\w\a\w dt+t^3 d\o_2\w\a-t^2d\o_3\w dt\\
=&-3t^2\o_2\w\a\w dt+t^3\nu_2\w\o_2\w\a-t^2(\nu_3\w\o_3+\sum_{r=1}^3
f_{3r}\o_r\w\a+\sigma_3\w\a)\w dt\,.
\end{aligned}
\]
Therefore $d\,\mbox{Re}\,\varepsilon=0$ if and only if one has the
extra-conditions
\[
f_{23}=3\,,\quad\nu_2=\nu_3=0\,,\quad \sigma_3=0\,.
\]
\end{proof}
\begin{rem}
\label{SI}
\em{
If $(\alpha,\omega_1,\omega_2,\omega_3)$ is a SU(2)-structure on $M$ inducing a symplectic half-flat
structure on $C(M)$, then defining
\[
\widetilde{\alpha}=\alpha\,,\quad\widetilde{\omega}_1=\omega_1\,,\quad\widetilde{\omega}_2=-\omega_3
\,,\quad\widetilde{\omega}_3=\omega_2\,,
\]
we obtain a {\em double-hypo} structure on $M$. The remarkable fact is that the two structures 
induce the same metric.
}
\end{rem}
\section{Ricci curvature of an $\mbox{SU}(2)$-structure}
Fix an SU(2)-reduction $\mathcal{Q}$ of the linear frame bundle
$\mathcal{L}(M)$, given by the quadruple
$(\alpha,\omega_1,\omega_2,\omega_3)$.  $\mathcal{Q}$ can be viewed
as a subbundle of the principal SO(5)-bundle $p\colon\mathcal{F}\to
M$ of the normal frames of the metric $g$ associated to the triple
$(\alpha,\omega_1,\omega_2,\omega_3)$. Consider on the bundle
$\mathcal{F}$ the tautological $\R^5$-valued $1$-form $w$ defined by
$w[u](v)=u(p_{*}[u]v)$ for every $u\in \mathcal{F}$ and $v\in
T_u\mathcal{F}$. On $\mathcal{F}$ we have also the Levi-Civita
connection 1-form $\psi$ taking values in $\mathfrak{so}(5)$. Using
the canonical basis $\{e_1,\dots,e_5\}$ of $\R^{5}$ we will regard
$w$ as a vector of $\R$-valued 1-forms on $\mathcal{F}$
$$
w=w_1e_1+\dots+w_5e_5
$$
and $\psi$ as a skew-symmetric matrix of 1-forms, i.e.
$\psi=(\psi_{ij})$. With this notation the first structure equation
relating $w$ and $\psi$
\begin{equation}
\label{1se} dw=-\psi\wedge w\,,
\end{equation}
becomes $d w_i=-\psi_{ij}\wedge w_j\,$. Note that equation
(\ref{1se})
simply means that $\psi$ is torsion-free.\\
The curvature of $\psi$ is by definition the
$\mathfrak{so}(5)$-valued 2-form $\Psi=d\psi+\psi\wedge\psi$. In
index notation
$$
\Psi_{ij}=d\psi_{ij}+\psi_{ik}\wedge\psi_{kj}=\frac{1}{2}\mathcal{R}_{ijkl}\,w_k\wedge
w_{l}\,.
$$
We consider the pull-backs of $\psi$ and $w$ to $\mathcal{Q}$  and
denote them by the same symbols for the sake of brevity. The
intrinsic torsion of the SU(2)-structure measures the failure of
$\psi$ to take values in $\mathfrak{su}$(2). More precisely
according to the splitting
$$
\mathfrak{so}(5)=
\mathfrak{su}(2)\oplus[\![\R^4]\!]\oplus[\R]_1\oplus[\R]_2\oplus[\R]_3\,.
$$
obtained above, $\psi$ decomposes as
$$
\psi=\theta+[\![\tau]\!]+[\mu_1]_1+[\mu_2]_2+[\mu_3]_3\,.
$$
Thus $\theta$ is a connection 1-form on $\mathcal{Q}$ which in
general is not torsion-free. We shall regard $\tau$ as a 4-vector of
1-forms $\tau=\tau_ie_i$. Furthermore we can write
$$
\tau_i=T_{ij}w_j\,,\quad \mu_r=M^r_jw_j
$$
for $i=1,2,3,4$ and $r=1,2,3$, where $T_{ij}$ and $M^r_i$ are smooth
functions. Formula \eqref{1se} now read as
$$
dw_i=-\theta_{ij}\wedge w_j-\eta_{ijk}\tau_k\wedge\
w_j-\epsilon^1_{ij}\,\mu_1\wedge w_j- \epsilon^2_{ij}\,\mu_2\wedge
w_j-\epsilon^3_{ij}\,\mu_3\wedge w_j\,,
$$
where the $\eta_{ijk}$'s are defined by \eqref{eta}.\\
Now we have
\begin{lemma}
\label{identita} The following identities hold:
\begin{enumerate}
\item[1.]$[\mu_r]_r\wedge[\![\tau]\!]+[\![\tau]\!]\wedge[\mu_r]_r=[\![[\mu_r]_r\wedge\tau]\!]\;$ for $r=1,2,3$;
\vskip0.3cm
\item[2.]$[\![\tau]\!]\wedge\theta+\theta\wedge[\![\tau]\!]=[\![\theta\wedge\tau]\!]$,
\end{enumerate}
where in the expressions $[\![[\mu_r]_r\wedge\tau]\!]$ and  $\theta\wedge\tau$, $\tau$ is regarded as 
the $\R^5$-valued $1$-form $\tau=(\tau_1,\dots,\tau_4,0)$.
\end{lemma}
We are ready to introduce the following quantities
$$
\begin{aligned}
& D\theta=d\theta+\theta\wedge\theta+[\![\tau]\!]\wedge[\![\tau]\!]
 +\frac{1}{4}\sum_{r=1}^3[\epsilon^{r}_{ij}\tau_i\wedge\tau_j]_r\,,\\
& D\tau=d\tau+\theta\wedge\tau+\sum_{r=1}^3[\mu_r]_r\wedge \tau\,,\\
& D\mu_1=d\mu_1-\frac{1}{4}\epsilon_{ij}^1\tau_i\wedge\tau_j-2\mu_2\wedge\mu_3\,,\\
& D\mu_2=d\mu_2-\frac{1}{4}\epsilon_{ij}^2\tau_i\wedge\tau_j-2\mu_3\wedge\mu_1\,,\\
&
D\mu_3=d\mu_3-\frac{1}{4}\epsilon_{ij}^3\tau_i\wedge\tau_j-2\mu_1\wedge\mu_2\,.
\end{aligned}
$$
A direct computation gives that $D\theta$ takes values in
$\mathfrak{su}(2)$; moreover lemma \ref{identita} implies
$$
\begin{aligned}
\Psi = & d(\theta+[\![\tau]\!]+[\mu_1]_1+[\mu_2]_2+[\mu_3]_3)\\
      & + (\theta+[\![\tau]\!]+[\mu_1]_1+[\mu_2]_2+[\mu_3]_3) \wedge
        (\theta+[\![\tau]\!]+[\mu_1]_1+[\mu_2]_2+[\mu_3]_3)\\
    = & D\theta+[\![D\tau]\!]+[D\mu_1]_1+[D\mu_2]_2+[D\mu_3]_3\,.
\end{aligned}
$$
In terms of the $w$-frame we shall write
$$
\begin{aligned}
&D\theta_{ij}=\frac{1}{2}S_{ijkl}w_k\wedge w_l\,,\\
&D\tau_i =\frac{1}{2}T_{ijk}w_j\wedge w_k\,,\\
&D\mu_r =\frac{1}{2}N^{r}_{kl}w_k\wedge w_l\,,
\end{aligned}
$$
where the coefficients are smooth functions such that
$$
\begin{aligned}
& S_{ijkl}=-S_{jikl}=-S_{ijlk}\,, \\
& T_{ijk}=-T_{ikj}\,,\\
& N^r_{kl}=-N^r_{lk}\,.
\end{aligned}
$$
In terms of the functions just introduced, the components of the
curvature tensor expresses as
$$
\mathcal{R}_{ijkl}=S_{ijkl}+\eta_{ijh}T_{hkl}+\epsilon^1_{ij}
N^1_{kl}+\epsilon^2_{ij}N^2_{kl}+\epsilon^3_{ij}N^3_{kl}\,,
$$
where the $\eta_{ijk}$'s are the symbols defined in \eqref{eta}.
Let Ric$_{ij}=\mathcal{R}_{ikkj}$ be the components of the Ricci
tensor of $g$. As an application of the Bianchi identities we have
the following theorem which gives a formula for the Ricci
tensor and the scalar curvature $s=\mbox{Ric}_{ii}$ of $g$ in terms of
intrinsic torsion.
\begin{theorem}
\label{senzaS}
The Ricci tensor does not depend on the functions $S_{ijkl}$ and
each component writes as
$$
\label{Ricciastorsion}
\mbox{\emph{Ric}}_{ij}=\sum_{r=1}^3\{\epsilon^r_{ik}N^r_{jk}+\epsilon^r_{jk}N^r_{ik}-\eta_{ijl}\epsilon^r_{lk}N^r_{k5}\}+
\delta_{i5}\delta_{j5}T_{kk5}+T_{ij5}\,.
$$
Consequently,
$$
\label{sastorsion}
s=2\sum_{r=1}^3(\epsilon^r_{ik}N^r_{ik})+2T_{kk5}\,.
$$
\end{theorem}
\subsubsection{The scalar curvature in terms of torsion forms}
Pulling back the structure forms to the SU(2)-bundle
$\pi\colon\mathcal{Q} \to M$, and using the frame $w_1,\ldots,w_5$,
one gets the standard expression for $\alpha, \omega_1, \omega_2,
\omega_3$:
\[
\pi^*(\alpha)=w_5\,,\quad\pi^*(\omega_r)=\frac{1}{2}\epsilon^r_{ij}w_i\wedge
w_j \mbox{ for } r=1,2,3\,.
\]

Applying the symmetries of the $\epsilon$-symbol, we have

\begin{prop}\label{derivatives}
The derivatives of the structure forms are
$$
\begin{aligned}
d\pi^*(\alpha)&=\tau_k\wedge w_k\,,\\
d\pi^*(\omega_1)&=\epsilon^1_{ij}\,\tau_i\wedge w_j\wedge
w_5-\epsilon^2_{ij}\,\mu_3\wedge w_i\wedge w_j+
\epsilon^3_{ij}\,\mu_2\wedge w_i\wedge w_j\,,\\
d\pi^*(\omega_2)&=\epsilon^2_{ij}\,\tau_i\wedge w_j\wedge
w_5-\epsilon^3_{ij}\,\mu_1\wedge w_i\wedge w_j+
\epsilon^1_{ij}\,\mu_3\wedge w_i\wedge w_j\,,\\
d\pi^*(\omega_3)&=\epsilon^3_{ij}\,\tau_i\wedge w_j\wedge
w_5-\epsilon^1_{ij}\,\mu_2\wedge w_i\wedge w_j+
\epsilon^2_{ij}\mu_1\wedge w_i\wedge w_j\,.
\end{aligned}
$$
\end{prop}
Proposition \ref{derivatives} allows to write down the pull-backs of
the torsion forms in terms of $T_{ij},\,M^{r}_i$. A direct
computation gives the following formulae
$$
\begin{aligned}
\pi^*(f_{11})&=\frac{1}{2}T_{ii}\,,\\
\pi^*(f_{12})&=\frac{1}{2}\epsilon^3_{ij}T_{ij}-2M^3_5\,,\\
\pi^*(f_{13})&=-\frac{1}{2}\epsilon^2_{ij}T_{ij}+2M^2_5\,,\\
\pi^*(f_{23})&=\frac{1}{2}\epsilon^1_{ij}T_{ij}-2M^1_5\,,\\
\pi^*(\phi_1)&=-\frac{1}{2}\epsilon^1_{ij}T_{ij}\,,\\
\pi^*(\phi_2)&=-\frac{1}{2}\epsilon^2_{ij}T_{ij}\,,\\
\pi^*(\phi_3)&=-\frac{1}{2}\epsilon^3_{ij}T_{ij}\,,\\
\pi^*(\nu_1)&=(2\epsilon^2_{ij}M^2_i+2\epsilon^3_{ij}M^3_i)\,w_j\,,\\
\pi^*(\nu_2)&=(2\epsilon^1_{ij}M^1_i+2\epsilon^3_{ij}M^3_i)\,w_j\,,\\
\pi^*(\nu_3)&=(2\epsilon^1_{ij}M^1_i+2\epsilon^2_{ij}M^2_i)\,w_j\,,\\
\pi^*(\nu_4)&=T_{i5}\,w_i\,,
\end{aligned}
$$
$$
\begin{aligned}
\pi^*(\sigma_1)&=\frac{1}{4}(\epsilon^1_{ip}(T_{pj}+T_{jp})+\epsilon^2_{ip}\epsilon^3_{qj}(T_{pq}+T_{qp}))\,w_i\wedge w_j\,,\\
\pi^*(\sigma_2)&=\frac{1}{4}(\epsilon^2_{ip}(T_{pj}+T_{jp})-\epsilon^1_{ip}\epsilon^3_{qj}(T_{pq}+T_{qp}))\,w_i\wedge w_j\,,\\
\pi^*(\sigma_3)&=\frac{1}{4}(\epsilon^3_{ip}(T_{pj}+T_{jp})+\epsilon^1_{ip}\epsilon^2_{qj}(T_{pq}+T_{qp}))\,w_i\wedge w_j\,,\\
\pi^*(\sigma_4)&=(T_{ji}+\frac{1}{2}\epsilon^r_{pq}\epsilon^r_{ij}T_{pq})\,w_i\wedge w_j+T_{i5}\,w_i \wedge w_5\,.
\end{aligned}
$$
\vskip0.3cm
\textbf{Warning:}
From now on we identify the structure and torsion forms with their pull-backs to the
principal SU(2)-bundle $\mathcal{Q}$.\\
\newline
Combining these formulae with \eqref{sastorsion} we get the following
\begin{theorem}
\label{scal}
The scalar curvature of the Riemannian metric induced by an
\emph{SU}$(2)$-structure with torsion
$(f_{ij},\phi_i,\nu_i,\sigma_i)$ on a $5$-manifold is
$$
\begin{aligned}
s=&-5f_{11}^2-\sum_{i=1}^3\phi_{i}^2-4\phi_{1}f_{23}+4\phi_{2}f_{13}-4\phi_{3}f_{12}+\sum_{i=1}^3d^*\nu_i-2d^*\nu_4
 -\sum_{i=1}^3\frac{1}{2}|\nu_i|^2\\&+\langle\nu_1,\nu_2\rangle
 +\langle\nu_1,\nu_3\rangle-\langle\nu_1,\nu_4\rangle+\langle\nu_2,\nu_3\rangle
-\langle\nu_2,\nu_4\rangle-\langle\nu_3,\nu_4\rangle-2*(df_{11}\wedge\omega_1^2)\\
&-\sum_{i=1}^4\frac{1}{2}|\sigma_i|^2\,.
\end{aligned}
$$
\end{theorem}
As a direct consequence of the previous theorem we have the
following characterization of the scalar curvature of some
special structures:
\begin{itemize}
\vskip0.1cm
\item Hypo manifolds: $s=-\phi_1^2-4\phi_1f_{23}-2|\nu_4|^2-\frac{1}{2}\sum_{i=2}^4|\sigma_i|^2$;
\vskip0.1cm
\item Contact-Hypo manifolds: $s=-4+8f_{23}-\frac{1}{2}|\sigma_2|^2-\frac{1}{2}|\sigma_3|^2$;
\vskip0.1cm
\item Double Hypo manifolds: $s=20-\frac{1}{2}|\sigma_2|^2-\frac{1}{2}|\sigma_4|^2$;
\vskip0.1cm
\item Sasaki-Einstein manifolds: $s=20$.
\end{itemize}
Hence we have
\begin{cor}
The scalar curvature of the metric induced by a double-Hypo
structure is always less or equal to $20$. Furthermore it is equal
to $20$ if and only if the double-Hypo structure is Sasaki-Einstein.
\end{cor}
\subsubsection{The Ricci curvature in terms of torsion forms}
According to the splitting \eqref{2-tensors} of symmetric 2-tensors,
the Ricci curvature of a metric $g$ associated to a SU(2)-structure
on a $5$-manifold decompose as follows
\begin{equation}
\label{splitRic}
\mbox{Ric}=\frac{\lambda}{4}\,g^T+\mu\,\alpha\otimes\alpha+\mbox{Ric}_0\,.
\end{equation}
We recall that the metric $g$ is said to be $\alpha$-\emph{Einstein} if
$$
\mbox{Ric}_0=0\,.
$$
(see e.g \cite{BGM}).

From the decomposition of the Ricci tensor \eqref{splitRic}, the  scalar curvature splits as
$$
s=\lambda+\mu\,.
$$
A straightforward computation gives the following formulae which express
$\lambda$ and $\mu$ in terms of torsion forms:
$$
\begin{aligned}
\lambda=&-4f_{11}^2-2\sum_{i=1}^3\phi_{i}^2-4\phi_{1}f_{23}+4\phi_{2}f_{13}-4\phi_{3}f_{12}+\sum_{i=1}^3d^*\nu_i-d^*\nu_4
 -\frac{1}{2}\sum_{i=1}^3|\nu_i|^2\\&+\langle\nu_1,\nu_2\rangle
 +\langle\nu_1,\nu_3\rangle-\langle\nu_1,\nu_4\rangle+\langle\nu_2,\nu_3\rangle
-\langle\nu_2,\nu_4\rangle-\langle\nu_3,\nu_4\rangle-|\sigma_4|^2\\[5pt]
&-*(df_{11}\wedge\omega_1^2)
\end{aligned}
$$
and
\begin{equation}\label{mu}
\begin{aligned}
\mu=&-f_{11}^2+\sum_{i=1}^3\phi_{i}^2-d^*\nu_4
 -\frac{1}{2}\sum_{i=1}^3|\sigma_i|^2+\frac{1}{2}|\sigma_4|^2-*(df_{11}\wedge\omega_1^2)\,.
\end{aligned}
\end{equation}
As a consequence of these formulae we get the following
\begin{prop}
Let $(M,\alpha,\omega_1,\omega_2,\omega_3)$ be a contact-Hypo
manifold. Assume that the Ricci tensor of  the metric induced by the
Hypo-structure satisfies
$$
\mbox{\emph{Ric}}(R_{\alpha},R_{\alpha})=4\,;
$$
then $(M,\alpha,\omega_1,\omega_2,\omega_3)$ is Sasaki $\alpha$-Einstein.
\end{prop}
\begin{proof}
For a Sasaki SU(2)-structure to be $\alpha$-Einstein is equivalent
to be Hypo (see Theorem 14 in \cite{CS}), so we only need to prove that
$(M,\alpha,\omega_1,\omega_2,\omega_3)$ is Sasaki.
By equations \eqref{CHtorsion}, in the contact-Hypo case formula \eqref{mu} reduces to
$$
\mu=4-\frac{1}{2}|\sigma_2|^2-\frac{1}{2}|\sigma_3|^2\,.
$$
Then condition $\mbox{Ric}(R_{\alpha},R_{\alpha})=4$ readily implies
$\sigma_2=\sigma_3=0$. Furthermore we have
$$
0=d^2\omega_2=df_{23}\wedge\omega_3\wedge\alpha
$$
which implies that $df_{23}=h\,\alpha$ for some $h\in
C^{\infty}(M,\R)$. Moreover
$$
0=d^2f_{23}=dh\wedge\alpha+h\,d\alpha=dh\wedge\alpha-2h\,\omega_1
$$
implies $h=0$. Hence $f_{23}$ is a constant function on $M$.\\
Let
$$
\widetilde{\varepsilon}=e^{(f_{23}-1)\log
t}\,(\omega_2+i\omega_3)\wedge(t\alpha+idt)\,,
$$
then $\widetilde{\varepsilon}$ is a
closed $(3,0)$-form on the almost K\"ahler cone $C(M)=M\times\R^+$.
Thus $C(M)$ is K\"ahler (see e.g. \cite{BV}, remark 1.1) and consequently
$(M,\alpha,\omega_1,\omega_2,\omega_3)$ is Sasaki.
\end{proof}
\begin{rem}{\em
Note that an SU(2)- structure satisfying the hypotheses of the proposition
above gives rise to an Einstein metric $g$ if and only if the scalar
curvature $s=\lambda+\mu$ is exactly 20.
Indeed, $g$ is Einstein if and only if $\mu=\frac{\lambda}{4}$ and the
hypothesis $\mbox{Ric}(R_{\alpha},R_{\alpha})=4$ means $\mu=4$.
}
\end{rem}
The main theorem is obtained using the following algorithm, analogous to the one used by Bryant in \cite{B2}:\\
\begin{itemize}
\item introduce the symbols $S_{ijk}$, $V^r_{ij}$ in the expressions  of the derivatives of the $T_{ij}$ and $M^r_i$:
\begin{eqnarray*}
&& dT_{ij}=T_{ik}\theta_{kj}+T_{kj}\theta_{ki}+S_{ijk}w_{k}\,,\\
&& dM_{i}^r=M^{r}_{k}\theta_{ki}+V^{r}_{ik}w_k\,.
\end{eqnarray*}
These symbols admit a geometric interpretation: for instance $S_{ijk}$'s keep track of the covariant derivative of
the $\mathfrak{so}(5)$-valued 1-form $[\![\tau]\!]$ with respect to the SU(2)-connection corresponding to the 1-form $\theta$:
\[
D_\theta [\![\tau]\!]=d[\![\tau]\!] + \theta \wedge [\![\tau]\!] + [\![\tau]\!] \wedge \theta\,.
\]
Analogously  $V^r_{ik}$'s keep track of the covariant derivative of
$[\mu_r]_r$ with respect to $\theta$.

\vskip0.3cm
\item write $T_{ijk}$ in terms of $T_{ij}$, $S_{ijk}$ and $M^r_i$; write $N^r_{ij}$ in terms of $M^r_i$, $T_{ij}$ and $V^r_j$.
This can be done since, for instance,
\[
\begin{aligned}
\qquad
D\tau_i = & \,d\tau_i+\theta_{ik}\wedge\tau_k+{\textstyle \sum _{r=1}^3\epsilon_{ik}^r\mu_r\w \tau_k}\\
        = & \,dT_{ij}\wedge w_j +T_{ij}d w_j+T_{kj}\theta_{ik}\wedge w_j+{\textstyle \sum_{r=1}^3\epsilon_{ik}^r\mu_r\w \tau_k}\\
        = & \,dT_{ij}\wedge w_j -T_{ik}\theta_{kj}\wedge w_j-T_{ik}\eta_{kjl}\tau_l\wedge w_j-
            {\textstyle \sum_{r=1}^3T_{ik}\epsilon^r_{kj}\mu_r\wedge w_j}\\
          & +T_{kj}\theta_{ik}\wedge w_j+{\textstyle \sum_{r=1}^3\epsilon_{ik}^r\mu_r\w \tau_k}\\
        = & \,dT_{ij}\wedge w_j -T_{ik}\theta_{kj}\wedge w_j-T_{kj}\theta_{ki}\wedge w_j-T_{ik}\eta_{kjl}\tau_l\wedge w_j\\
          & -{\textstyle \sum_{r=1}^3(T_{ik}\epsilon^r_{kj}\mu_r \wedge w_j+\epsilon_{ik}^r\mu_r\w \tau_k)}\\
        = & \,S_{ijk}w_k\w w_j-T_{ik}\eta_{kjl}\tau_l\wedge w_j
            -{\textstyle \sum_{r=1}^3(T_{ik}\epsilon^r_{kj}\mu_r \wedge w_j+\epsilon_{ik}^r\mu_r\w \tau_k)}\\
         = & \,S_{ijk}w_k\w w_j-T_{im}T_{lk}\eta_{mjl}w_k\wedge w_j\\
           & -{\textstyle \sum_{r=1}^3(T_{il}M^r_k\epsilon^r_{lj}w_k \wedge w_j+T_{lj}M^r_k\epsilon_{il}^r w_k \w w_j)}\,,\\
\end{aligned}         
\]
i.e.
\[
\quad \quad \quad
T_{ikj}w_k \w w_j= \,(S_{ijk}-T_{im}T_{lk}\eta_{mjl}
           -{\textstyle \sum_{r=1}^3(T_{il}M^r_k\epsilon^r_{lj}+T_{lj}M^r_k\epsilon_{il}^r)})w_k \w w_j\,.\\
\]
\vskip0.3cm
\item  use Theorem \eqref{Ricciastorsion} to write the tensor
$\mbox{Ric}_0$ in terms of $T_{ij}, S_{ijk}, M^r_{i},V^r_{ij}$. The resulting expression
is linear in $S_{ijk},V^{r}_{ij}$ and  at most quadratic in $T_{ijk},M^r_{i}$;
\vskip0.3cm
\item decompose Ric$_0$ in
$$
\mbox{Ric}_0=\mbox{Ric}_0^{(1)}+\mbox{Ric}_0^{(2)}+\mbox{Ric}_0^{(3)}+\Phi_4\odot\alpha\,.
$$
according to splitting \eqref{decsym} and use the isomorphisms $\iota_r$'s to make Ric$_0$ into a  2-form
$\Phi=\Phi^{(1)}+\Phi^{(2)}+\Phi^{(3)}+\Phi_4\wedge\alpha$;
\vskip0.3cm
\item use representation theory of SU(2) to build the expressions, bilinear in $\{\nu_i,\sigma_j,\phi_r,f_{uv}\}$
and linear in their derivatives, sufficient to write $\Phi$ as a linear combination of them.
\end{itemize}

\begin{theorem}\label{Ricci}
The ``traceless part'' of the Ricci curvature of the Riemannian
metric induced by an \emph{SU}$(2)$-structure with torsion
$(f_{ij},\phi_i,\nu_i,\sigma_i)$ on a $5$-manifold is
$$
\mbox{\emph{Ric}}_0=\iota_1^{-1}(\mbox{\emph{E}}(\Phi_{1}))+\iota_2^{-1}(\mbox{\emph{E}}(\Phi_{2}))
+\iota_3^{-1}(\mbox{\emph{E}}(\Phi_{3}))+\Phi_4\odot\alpha\,,
$$
where
\[
\begin{array}{rcl}
\Phi_{1}&=&-\frac{1}{2}f_{11}\sigma_1+\frac{1}{2}f_{12}\sigma_2+\frac{1}{2}f_{13}\sigma_3
-f_{23}\sigma_4+\phi_3\sigma_2-\phi_2\sigma_3-\phi_{1}\sigma_4-\frac{1}{4}
\nu_1\wedge J_1\nu_1\\[4pt]
& &+\frac{1}{2}\nu_1\wedge J_1\nu_4+\frac{1}{4}\nu_2\wedge
J_1\nu_2-\frac{1}{2}\nu_2\wedge J_1\nu_3+\frac{1}{4}\nu_3\wedge
J_1\nu_3-\frac{1}{2}\nu_2\wedge J_1\nu_3\\[4pt]
& &+\frac{1}{2}\nu_4\wedge J_1\nu_4
+\frac{1}{2}\iota_{R_{\alpha}}d\sigma_1
-\frac{1}{2}dJ_1\nu_1+\frac{1}{2}dJ_1\nu_2+\frac{1}{2}dJ_1\nu_3\,;\\[10pt]
\Phi_{2}&=&-\frac{1}{2}f_{12}\sigma_2-\frac{1}{2}f_{11}\sigma_2+\frac{1}{2}f_{23}\sigma_3
+f_{13}\sigma_4-\phi_3\sigma_1 +\phi_1\sigma_3-\phi_{2}\sigma_4
+\frac{1}{2}*d\sigma_2\\[4pt]
& &+\frac{1}{2}dJ_2\nu_1
-\frac{1}{2}dJ_2\nu_2-\frac{1}{2}dJ_2\nu_4+\frac{1}{2}dJ_2\nu_3
+\frac{1}{4}\nu_1\wedge J_2\nu_1-\frac{1}{2}\nu_1\wedge
J_2\nu_3\\[4pt]
& &-\frac{1}{4}\nu_2\wedge J_2\nu_2+\frac{1}{2}\nu_2\wedge
J_2\nu_4+\frac{1}{2}\nu_4\wedge J_2\nu_4+\frac{1}{4}\nu_3\wedge
J_2\nu_3\,;\\[10pt]
\Phi_3 & = &
-\frac{1}{2}(f_{13}\sigma_1+f_{23}\sigma_2+f_{11}\sigma_3)+\phi_2\sigma_1-\phi_1\sigma_2-\phi_3\sigma_4
+\frac{1}{2}*d\sigma_3-f_{12}\sigma_4\\[4pt]
& &+\frac{1}{4}(\nu_1\wedge J_3\nu_1+\nu_2\wedge
J_3\nu_2-\nu_3\wedge J_3\nu_3)-\frac{1}{2}\nu_1\wedge J_3
\nu_2+\frac{1}{2}\nu_3\wedge J_3\nu_4\\[4pt]
& &+\frac{1}{2}\nu_4\wedge J_3\nu_4+\frac{1}{2}(dJ_3\nu_1+dJ_3\nu_2-dJ_3\nu_3-dJ_3\nu_4)\,;\\[10pt]
\Phi_4 & = &
3(df_{11})^T-\frac{3}{2}f_{11}\nu_4-\frac{1}{2}(d^*\sigma_4)^T-\frac{1}{2}J_2(d^*\sigma_2)^T-\frac{1}{2}f_{23}J_1\nu_4
-\frac{1}{2}f_{12}J_3\nu_4\\[4pt]
& & +\frac{3}{2}(\phi_1J_1\nu_4+\phi_2 J_2\nu_4+\phi_3 J_3\nu_4)-\frac{1}{2}\iota_{R_\alpha}(d\nu_1+d\nu_2+d\nu_3)\\[4pt]
& &
+\iota_{R_{\alpha}}(d\nu_4+*(df_{12}\wedge\omega_3)-*(df_{13}\wedge\omega_2)+*(df_{23}\wedge\omega_1)+
\frac{1}{2}*d\sigma_4)\\[4pt]
& &-J_1\iota_{R_{\alpha}}(*d\sigma_1)
-\frac{3}{2}J_2\iota_{R_{\alpha}}(*d\sigma_2)-J_{3}\iota_{R_{\alpha}}(*d\sigma_3)
+\frac{1}{2}J_1\iota_{R_{\alpha}}(dJ_1\nu_4)\\[4pt]
& &+\frac{1}{2}J_3\iota_{R_{\alpha}}(dJ_3\nu_4)\,;\\[4pt]
\end{array}
\]
and the operators $\iota_r\colon\Sigma_r(M) \to \Lambda^2_3M$ and
$\mbox{\emph{E}}\colon\Lambda^2M \to \Lambda^2_3M$ are
defined respectively in \eqref{iota} and \eqref{E}.
\end{theorem}
\section{The Ricci tensor of a contact-Hypo manifold}
\noindent Let $(M,\alpha,\omega_1,\omega_2,\omega_3)$ be a
contact-Hypo manifold. In view of the observations of subsection \ref{akcone}, $\Phi_1,\Phi_2,\Phi_3,\Phi_4$
reduce to
\[
\begin{aligned}
\Phi_1=&0\,,\\
\Phi_2=&(\frac{1}{2}f_{23}-2)\sigma_3+\frac{1}{2}*d\sigma_2\,,\\
\Phi_3=&(-\frac{1}{2}f_{23}+2)\sigma_2+\frac{1}{2}*d\sigma_3\,,\\
\Phi_4=&-\frac{1}{2}J_2(d^*\sigma_2)^T+\iota_{R_{\alpha}}(*(df_{23}\wedge\omega_1))-\frac{3}{2}J_{2}\iota_{R_{\alpha}}(*d\sigma_2)-J_{3}
\iota_{R_{\alpha}}(*d\sigma_3)\,.
\end{aligned}
\]
Now we observe that
\begin{eqnarray}
&& \mbox{E}((\frac{1}{2}f_{23}-2)\sigma_3+\frac{1}{2}*d\sigma_2)=(\frac{1}{2}f_{23}-2)\sigma_3+\frac{1}{2}\mbox{E}(*d\sigma_2)\,,\\
&&
\mbox{E}((-\frac{1}{2}f_{23}+2)\sigma_2+\frac{1}{2}*d\sigma_3)=(-\frac{1}{2}f_{23}+2)\sigma_2+\frac{1}{2}\mbox{E}(*d\sigma_3)\,.
\end{eqnarray}
Moreover, using \eqref{prop*}, for $i=2,3$,  we get
\begin{eqnarray*}
\mbox{E}(*d\sigma_i)=*^T\mbox{E}(\iota_{R_{\alpha}}d\sigma_i)=*^T(\iota_{R_{\alpha}}d\sigma_i)-\frac{1}{2}\sum_{r=1}^3*(\iota_{R_{\alpha}}d\sigma_i
\wedge\omega_r\wedge\alpha)\,\omega_r\,.
\end{eqnarray*}
Consequently
\begin{eqnarray}
\label{EP1EP2} && \mbox{E}(\Phi_2)=(\frac{1}{2}f_{23}-2)\sigma_3+
\frac{1}{2}*^T(\iota_{R_{\alpha}}d\sigma_2)-\frac{1}{4}\sum_{r=1}^3*(\iota_{R_{\alpha}}d\sigma_2\wedge\omega_r\wedge\alpha)\,\omega_r\,,\\
&& \mbox{E}(\Phi_3)=
(-\frac{1}{2}f_{23}+2)\sigma_2+\frac{1}{2}*^T(\iota_{R_{\alpha}}d\sigma_3)-
\frac{1}{4}\sum_{r=1}^3*(\iota_{R_{\alpha}}d\sigma_3\wedge\omega_r\wedge\alpha)\,\omega_r\,.
\end{eqnarray}
In order to write down the Ricci tensor of a contact-Hypo structure,
we consider the following
\begin{lemma}
Let $(M,\alpha,\omega_1,\omega_2,\omega_3)$ be a contact-Hypo
manifold, then
$$
\Phi_4=3 J_1(df_{23})^T\,.
$$
\end{lemma}
\begin{proof}
The lemma is essentially a consequence of the vanishing of
$d^{\,2}$. First of all note that, for any 3-form $\gamma$, one can
write $-*^T\gamma^T$ instead of $\iota_{R_\alpha}*\gamma$. Hence in
the contact-Hypo case one has
\begin{equation}
\label{Fi4}
\Phi_4=-\frac{1}{2}J_2(d^*\sigma_2)^T-*^T(df_{23}\wedge\omega_1)^T+\frac{3}{2}J_2*^T(d\sigma_2)^T+J_{3}*^T(d\sigma_3)^T\,.
\end{equation}
Now
\[
\begin{aligned}
&0=d^2\omega_2=-\alpha\wedge(df_{23}\wedge\omega_3+d\sigma_2)\,,\\
&0=d^2\omega_3=-\alpha\wedge(-df_{23}\wedge\omega_2+d\sigma_3)\,.
\end{aligned}
\]
Hence
\[
\begin{aligned}
&(d\sigma_2)^T=-(df_{23}\wedge \omega_3)^T=-(df_{23})^T\omega_3=-J_3*^T(df_{23})^T\,,\\
&(d\sigma_3)^T=(df_{23}\wedge
\omega_2)^T=(df_{23})^T\omega_2=J_2*^T(df_{23})^T\,.
\end{aligned}
\]
For the first term of \eqref{Fi4}, consider
\[
d*\sigma_2=-d(\sigma_2\wedge\alpha)=-d\sigma_2\wedge\alpha+2\,\sigma_2\wedge\omega_1=-(d\sigma_2)^T\wedge\alpha\,.
\]
Thus
\[
J_2(d^*\sigma_2)^T=J_2*(\alpha\wedge(d\sigma_2)^T)=J_2*^T(d\sigma_2)^T=-J_2J_3(*^T)^2(df_{23})^T=J_1(df_{23})^T\,.
\]
Finally for the second term
\[
*^T(df_{23}\wedge\omega_1)^T=*^T((df_{23})^T\wedge\omega_1)=*^TJ_1*^T(df_{23})^T=-J_1(df_{23})^T\,.
\]
Therefore, keeping in mind the quaternionic relations of $J_r$'s,
one has
\[
\Phi_4=(-\frac{1}{2}+1+\frac{3}{2}+1)J_1(df_{23})^T=3J_1(df_{23})^T\,.
\]
\end{proof}
Summarizing we have the following
\begin{prop}
The ``traceless part'' of the Ricci tensor of a
contact-Hypo manifold is given by the following formula
\begin{equation}
\label{riccihypo}
\begin{aligned}
{\rm Ric}_0=&\,\iota_2^{-1}\Big((\frac{1}{2}f_{23}-2)\sigma_3+
\frac{1}{2}*^T(\iota_{R_{\alpha}}d\sigma_2)-\frac{1}{4}\sum_{r=1}^3*(\iota_{R_{\alpha}}d\sigma_2\wedge\omega_r\wedge\alpha)\,\omega_r\Big)\\
&+\iota_3^{-1}\Big((-\frac{1}{2}f_{23}+2)\sigma_2+\frac{1}{2}*^T(\iota_{R_{\alpha}}d\sigma_3)-
\frac{1}{4}\sum_{r=1}^3*(\iota_{R_{\alpha}}d\sigma_3\wedge\omega_r\wedge\alpha)\,\omega_r\Big)\\
&+ 3 J_1(df_{23})^T\odot\alpha\,.
\end{aligned}
\end{equation}
\end{prop}
Now we collect some consequences of this result.\\

\begin{prop}
\label{double}
Let $(M,\alpha,\omega_1,\omega_2,\omega_3)$ be a double-Hypo
$5$-manifold. The metric induced by the
$\mbox{\emph{SU}}(2)$-structure is $\alpha$-Einstein if and only if
$(M,\alpha,\omega_1,\omega_2,\omega_3)$ is Sasaki-Einstein.
\end{prop}
\begin{proof}
The $\alpha$-Einstein condition means that the projection onto $\Lambda^2_3 M$ of 
$\Phi_1, \Phi_2$ and $\Phi_3$ vanishes.
But in the double-Hypo case one has
\[
\Phi_1=-\sigma_4\,, \quad \Phi_2=-\frac{1}{2}\sigma_3\,,
\]
which lie in  $\Lambda^2_3 M$, and the conclusion follows.
\end{proof}

\begin{cor}
\label{SGCY}
Let $(M,\alpha,\omega_1,\omega_2,\omega_3)$ be a $5$-dimensional
\emph{SU}$(2)$-manifold.\\ Assume that:
\begin{enumerate}
\item[1.] the \emph{SU}$(3)$-structure induced on the cone $C(M)=M\times\R^+$ is symplectic half-flat,
\vskip0.2cm
\item[2.] the metric $g$ induced by $(\alpha,\omega_1,\omega_2,\omega_3)$ is $\alpha$-Einstein,
\end{enumerate}
then $(M,\alpha,\omega_1,\omega_2,\omega_3)$ is Sasaki-Einstein.
\end{cor}
\begin{proof}
Simply recall remark \ref{SI} and apply the previous proposition.  
\end{proof}

\end{document}